\documentclass[reqno,12pt,letterpaper]{amsart}
\usepackage{sdmacros}

\DeclareMathOperator{\Sp}{Sp}

\begin{document}

\author{Semyon Dyatlov}
\email{dyatlov@math.mit.edu}
\address{Department of Mathematics, Massachusetts Institute of Technology, Cambridge, MA 02139}
\title{Macroscopic limits of chaotic eigenfunctions}

\begin{abstract}
We give an overview of the interplay between the behavior
of high energy eigenfunctions of the Laplacian on a compact Riemannian manifold and
the dynamical properties of the geodesic flow on that manifold. This includes the Quantum Ergodicity theorem,
the Quantum Unique Ergodicity conjecture, entropy bounds, and uniform lower bounds on mass of eigenfunctions.
The above results belong to the domain of \emph{quantum chaos} and use \emph{microlocal analysis}, which is a theory behind the classical/quantum, or particle/wave, correspondence in physics. We also discuss the toy model of quantum cat maps
and the challenges it poses for Quantum Unique Ergodicity.
\end{abstract}

\maketitle

\section{Introduction}

This article is an overview of some results on \emph{macroscopic behavior
of eigenstates in the high energy limit}.
A typical model is given by Laplacian eigenfunctions:
\[
-\Delta_g u_\lambda=\lambda^2u_\lambda,\qquad
u_\lambda\in C^\infty(M),\qquad
\lVert u_\lambda\rVert_{L^2(M)}=1.
\]
Here we fix a compact connected Riemannian manifold without boundary $(M,g)$
and denote by $\Delta_g\leq 0$ the corresponding Laplace--Beltrami operator.
It will be convenient to denote the eigenvalue by $\lambda^2$, where
$\lambda\geq 0$. The high energy limit corresponds to taking $\lambda\to\infty$.

One way to study macroscopic behavior of the eigenfunctions $u_\lambda$
as $\lambda\to\infty$ is to look at weak limits of the probability measures
$|u_\lambda|^2\,d\vol_g$ where $d\vol_g$ is the volume measure on~$(M,g)$:
\begin{defi}
\label{d:weak-limit-1}
Let $\lambda_j^2$ be a sequence of eigenvalues of $-\Delta_g$ going to~$\infty$.
We say that the corresponding eigenfunctions $u_{\lambda_j}$ converge weakly
to some probability measure $\nu$ on~$M$, if
\begin{equation}
  \label{e:weak-limit-1}
\int_M a(x)|u_{\lambda_j}(x)|^2\,d\vol_g(x)\ \to\ \int_M a(x)\,d\nu(x)\quad\text{as}\quad
j\to\infty
\end{equation}
for all test functions $a\in C^\infty(M)$.
\end{defi}
Definition~\ref{d:weak-limit-1} can be interpreted in the context of quantum mechanics
as follows. Consider a free quantum particle on the manifold $M$. Then the eigenfunctions $u_\lambda$ are the wave functions of the \emph{pure quantum states} of the particle.
The left-hand side of~\eqref{e:weak-limit-1} is the average value of the observable
$a(x)$ for a given pure state; if we let $a$ be the characteristic function
of some set $\Omega\subset M$ then this expression is the probability of finding
the quantum particle in $\Omega$ (this choice is only allowed if $\nu(\partial\Omega)=0$).
Taking $\lambda\to\infty$ gives the high energy limit.

The statement~\eqref{e:weak-limit-1}
is macroscopic in nature because we first fix the observable~$a$ and then let
the eigenvalue go to infinity. This is different from \emph{microscopic} properties
such as the breakthrough work of Logunov and Malinnikova on the area of the \emph{nodal set}
$\{x\in M\mid u_j(x)=0\}$, see the review~\cite{Logunov-Malinnikova-Yau-review}.
Ironically the methods used in the macroscopic results described here are \emph{microlocal} in nature
(see~\S\ref{s:semi-measures} for a review), with the global geometry of~$M$ coming in the form
of the long time behavior of the geodesic flow.

The results reviewed in this paper address the following fundamental question:
\begin{equation}
\label{e:question-1}
\begin{gathered}
\text{For a given Riemannian manifold $(M,g)$, what can we say}\\
\text{about the set of all weak limits of sequences of eigenfunctions?}
\end{gathered}
\end{equation}
\begin{figure}
\hbox to\hsize{
\hss
\includegraphics[width=6cm]{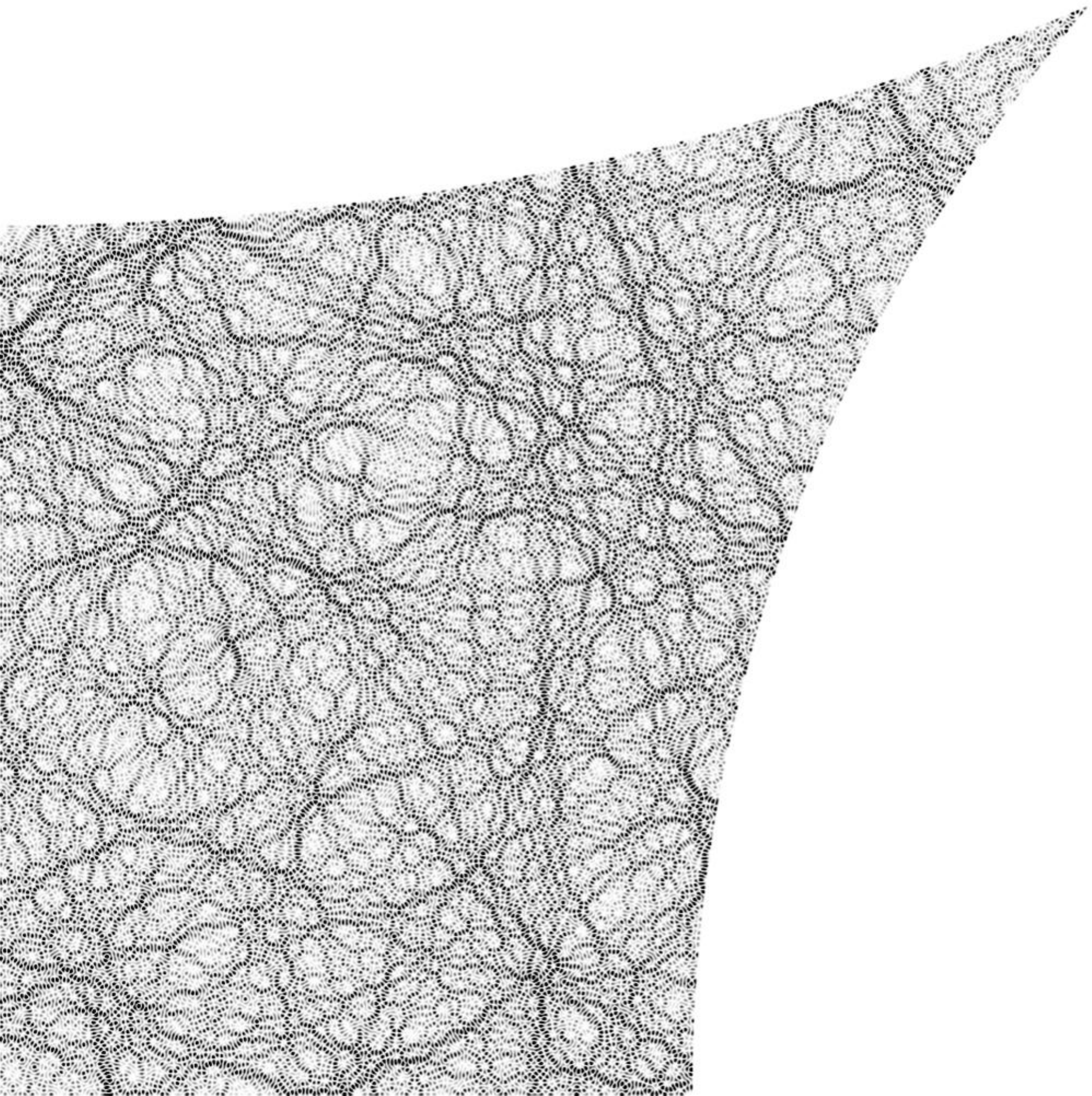}\qquad
\includegraphics[width=6cm]{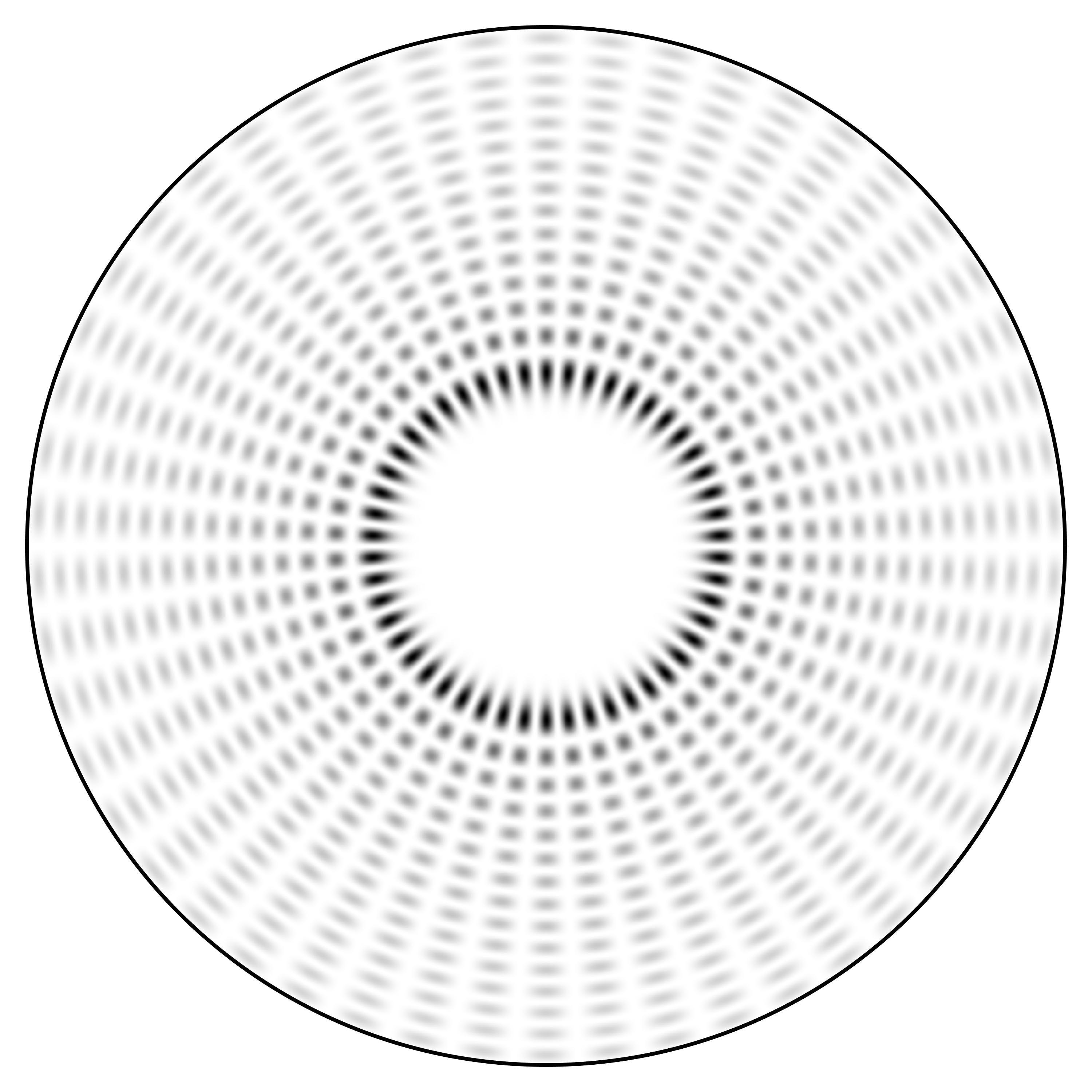}
\hss}
\hbox to\hsize{
\hss
\includegraphics[width=6cm]{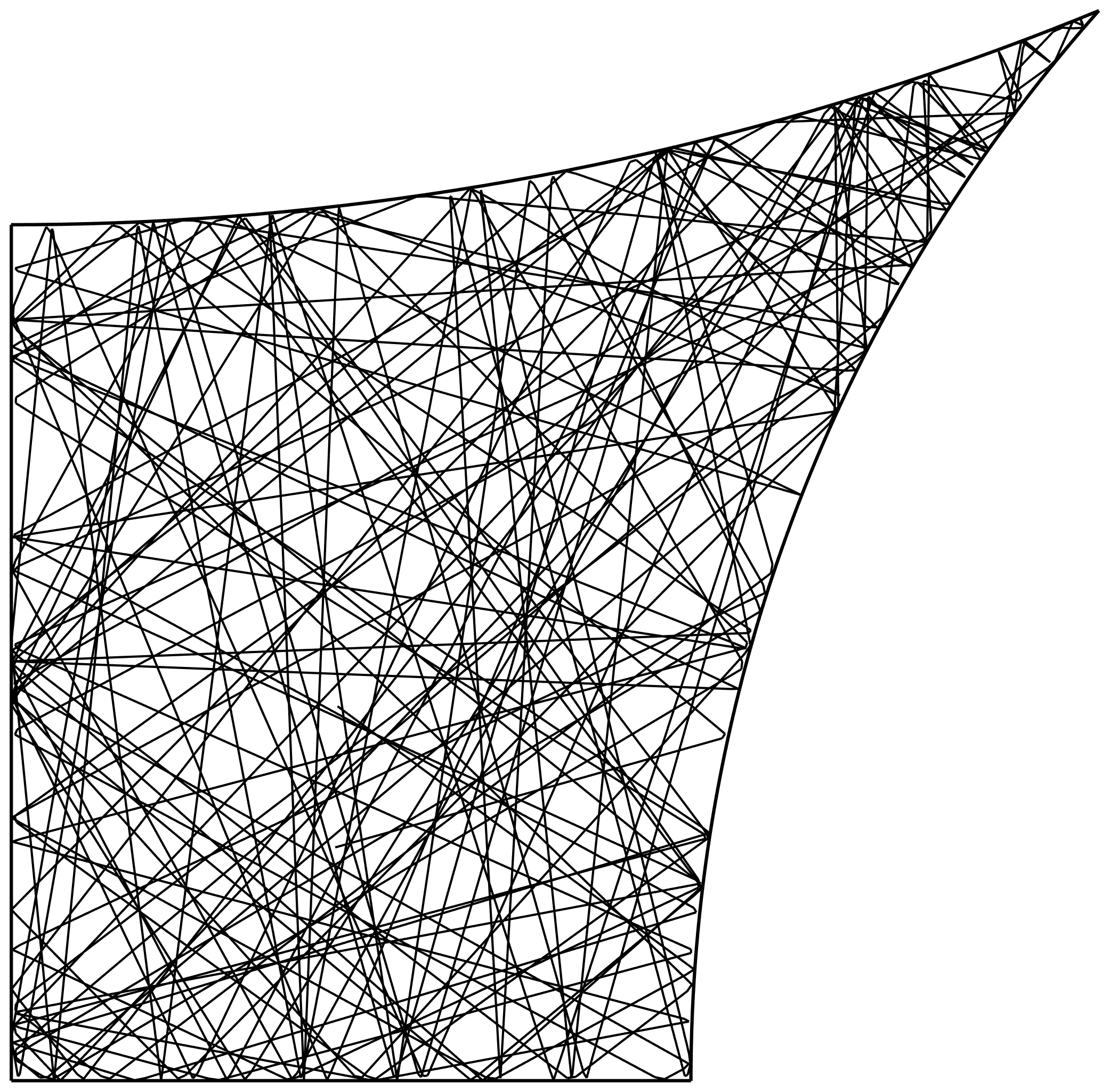}\qquad
\includegraphics[width=6cm]{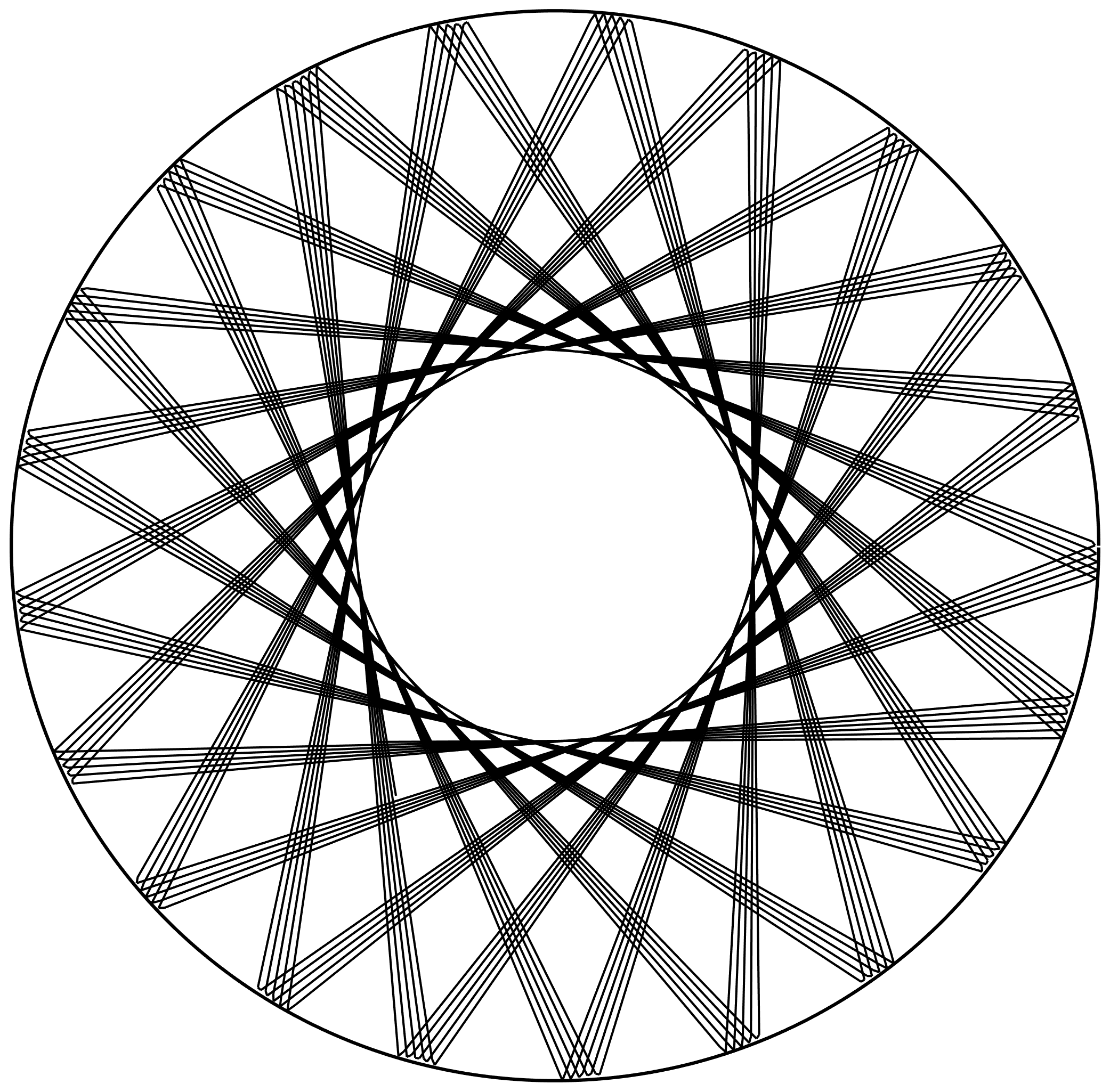}
\hss}
\caption{Top: typical eigenfunctions (with Dirichlet boundary conditions)
for two planar domains. The picture on the left (courtesy of Alex Barnett,
see~\cite{Barnett-Billiard} and~\cite{Barnett-Hassell} for a description of the method
used and for a numerical investigation of Quantum Ergodicity, showing empirically
$\mathcal O(\lambda^{-1/2})$ convergence to equidistribution)
shows equidistribution, i.e. convergence to the volume measure in the sense of Definition~\ref{d:weak-limit-1}. The picture on the right (where the domain is a disk)
shows lack of equidistribution, with the limiting measure supported in an annulus.
This difference in quantum behavior is related to the different behavior
of the billiard-ball flows on the two domains (which replace geodesic flows
in this setting). Bottom: two typical billiard-ball trajectories
on the domains in question. On the left we see ergodicity (equidistribution of the
trajectory for long time) and on the right we see completely integrable behavior.
}
\label{f:billiards}
\end{figure}
It turns out that the answer depends on the dynamical properties of the \emph{geodesic flow}
on~$(M,g)$. In particular:
\begin{itemize}
\item If $(M,g)$ has \emph{completely integrable} geodesic flow
then there is a huge variety of possible weak limits. For example,
if $(M,g)$ is the round sphere, then there is a sequence of Gaussian
beam eigenfunctions converging to the delta measure on any given closed geodesic
(see~\S\ref{s:semimes} below).
\item If the geodesic flow instead has \emph{chaotic} behavior, more precisely it is ergodic
with respect to the Liouville measure, then a density one sequence of eigenfunctions
converges to the volume measure $d\vol_g/\vol_g(M)$. This statement,
known as \emph{Quantum Ergodicity}, is reviewed in~\S\ref{s:QE}.
\item If the geodesic flow is \emph{strongly chaotic}, more precisely it satisfies the Anosov property (i.e. it has a stable/unstable/flow decomposition), then the limiting measures
have to be somewhat spread out. This comes in two forms: \emph{entropy bounds}
and \emph{full support}. See~\S\ref{s:anosov} for a description of these results.
The \emph{Quantum Unique Ergodicity} conjecture states that in this setting
any sequence of eigenfunctions converges to the volume measure;
it is not known outside of arithmetic cases (see~\S\ref{s:anosov}) and there are
counterexamples in the related setting of quantum cat maps (see~\S\ref{s:quantum-cat}).
\item Finally, there are several results in cases when the geodesic flow is ergodic but not Anosov, or it exhibits mixed chaotic/completely integrable behavior~-- see~\S\ref{s:QE}.
\end{itemize}
The present article focuses on the last three cases above, which are in the domain
of \emph{quantum chaos}. The general principle is that \emph{chaotic behavior of the geodesic
flow leads to chaotic/spread out macroscopic behavior of the eigenfunctions of the Laplacian}. See Figure~\ref{f:billiards} for a numerical illustration.

In particular, we will describe full support statements for weak limits~-- see Theorem~\ref{t:meassupp}
and Theorem~\ref{t:highcat}~-- proved in~\cite{meassupp,varfup,highcat}. The key component is the
\emph{fractal uncertainty principle} first introduced by Dyatlov--Zahl~\cite{hgap} and
proved by Bourgain--Dyatlov~\cite{fullgap}. It originated in \emph{open} quantum chaos, dealing
with quantum systems where the underlying classical system allows escape to infinity
and has chaotic behavior. We refer to the reviews of the author~\cite{Dyatlov-JEDP,FUP-ICMP} for more on fractal uncertainty principle and its applications.

The above developments use \emph{microlocal analysis}, which is a mathematical theory
underlying the classical/quantum, or particle/wave, correspondence in physics. In particular,
one typically obtains information on the \emph{semiclassical measures},
which are probability measures $\mu$ on the cosphere bundle $S^*M$ which are weak limits of
sequences of eigenfunctions in a microlocal sense. These measures are sometimes called
\emph{microlocal lifts} of the weak limits, because the pushforward
of $\mu$ to the base~$M$ is the weak limit of Definition~\ref{d:weak-limit-1}.
One of the advantages of these measures compared to the weak limits on~$M$ is that they
are invariant under the geodesic flow.
We give a brief review of microlocal analysis and semiclassical measures
in~\S\ref{s:semi-measures} below.

\section{Semiclassical measures}
  \label{s:semi-measures}
  
Let us write the left-hand side of~\eqref{e:weak-limit-1} as
\[
\int_M a(x)|u_{\lambda_j}(x)|^2\,d\vol_g(x)
=\langle \mathbf M_a u_{\lambda_j},u_{\lambda_j}\rangle_{L^2(M)}
\]
where $\mathbf M_a:L^2(M)\to L^2(M)$ is the multiplication operator
by $a\in C^\infty(M)$. To define semiclassical measures we will allow
more general operators in place of $\mathbf M_a$. These operators
are obtained by a \emph{quantization procedure},
which maps each smooth compactly supported function $a$ on the cotangent bundle~$T^*M$
to an operator on $L^2(M)$ depending on the small number~$h>0$
called the semiclassical parameter:
\begin{equation}
  \label{e:quant-proc}
a\in \CIc(T^*M)\quad \mapsto\quad \Op_h(a):L^2(M)\to L^2(M),\quad
0<h\ll 1.
\end{equation}

\subsection{Semiclassical quantization}

We briefly recall several basic principles of semiclassical quantization referring to
the books of Zworski~\cite{Zworski-Book}
and Dyatlov--Zworski~\cite[Appendix~E]{DZ-Book} for full presentation and pointers to the vast literature on the subject:
\begin{itemize}
\item The function $a$, often called the \emph{symbol} of the operator $\Op_h(a)$,
is defined on the cotangent bundle $T^*M$, whose
points we typically denote by $(x,\xi)$ where $x\in M$ and $\xi\in T_x^*M$.
The canonical symplectic form on $T^*M$
induces the \emph{Poisson bracket}
\[
\{f,g\}:=\partial_\xi f\cdot\partial_xg-\partial_x f\cdot\partial_\xi g, \quad
f,g\in C^\infty(T^*M).
\]
In physical terms, this corresponds to using
Hamiltonian mechanics for the `classical' side of the classical/quantum correspondence,
where $x$ is the position variable and $\xi$ is the momentum variable.
\item One can work with a broader class of smooth
symbols $a$, where the compact support requirement is changed to growth conditions
on the derivatives of~$a$ as $\xi\to\infty$. The resulting
operators act on (semiclassical) Sobolev spaces, see e.g.~\cite[\S E.1.8]{DZ-Book}. 
\item If $a(x,\xi)=a(x)$ is a function of $x$ only, then
\begin{equation}
  \label{e:op-h-mul}
\Op_h(a)=\mathbf M_a
\end{equation}
is the corresponding multiplication operator.
\item If $a(x,\xi)$ is linear in~$\xi$, that is $a(x,\xi)=\langle \xi,X_x\rangle$
for some vector field $X\in C^\infty(M;TM)$, then up to lower order terms
the operator $\Op_h(a)$ is a rescaled differentiation operator along $X$:
\begin{equation}
  \label{e:op-h-vf}
\Op_h(a) u(x)=-ihXu(x)+\mathcal O(h).
\end{equation}
This explains why $a$ should be a function on the cotangent bundle $T^*M$:
linear functions on the fibers of $T^*M$ correspond to vector fields on~$M$.
(Quantization procedures do not depend on the choice of a Riemannian metric on~$M$.)
\item If $u\in C^\infty(M)$
oscillates at some frequency $R$, then differentiating $u$ along a vector field $X$
increases its magnitude by about $R$. 
One takeaway from~\eqref{e:op-h-vf} is that $\Op_h(a)u$ has roughly the
same size as $u$ if the function $u$ oscillates at frequencies $\sim h^{-1}$.
Thus we treat the semiclassical parameter $h$ as the \emph{effective wavelength}
of oscillations of the functions to which we will apply $\Op_h(a)$.
We will apply $\Op_h(a)$ to an eigenfunction $u_\lambda$, which oscillates
at frequency $\sim \lambda$, so we will make the choice
\begin{equation}
\label{e:h-chosen}
h:=\lambda^{-1}.
\end{equation}
\item If $M=\mathbb R^n$ and $a(x,\xi)=a(\xi)$ is a function of~$\xi$ only, then
$\Op_h(a)$ is a Fourier multiplier:
\begin{equation}
  \label{e:op-h-fourier}
\widehat{\Op_h(a) u}(\xi)=a(h\xi)\hat u(\xi),\quad
u\in\mathscr S(\mathbb R^n).
\end{equation}
Thus in addition to being the momentum variable, we can interpret $\xi$
as a Fourier/frequency variable.
\item For general manifolds $M$, one cannot define a quantization procedure canonically:
a typical construction involves piecing together quantizations on copies of~$\mathbb R^n$
using coordinate charts, see e.g.~\cite[\S E.1.7]{DZ-Book}. However,
different choices of coordinate charts etc. will give the same operator
modulo lower order terms $\mathcal O(h)$.
\end{itemize}
Several items above allude to `lower order terms'. We will consider the operators
$\Op_h(a)$ in the \emph{semiclassical limit} $h\to 0$ and will often have
remainders of the form $\mathcal O(h)$ etc. which are operators on $C^\infty(M)$.
(More generally, semiclassical analysis gives asymptotic expansions in powers
of $h$ with remainder being $\mathcal O(h^N)$ for any~$N$.) This is understood
as follows: if the symbols
involved are compactly supported in $T^*M$, then the remainders
are bounded in norm as operators on $L^2$ (with constants in $\mathcal O(\bullet)$
of course independent of~$h$). For more general symbols, one has to take
correct semiclassical Sobolev spaces and we skip these details here.
We note that in the basic version of semiclassical calculus used
in this section, the symbol~$a$ does not depend on~$h$, which
reflects the macroscopic nature of the results presented below.

Semiclassical quantization has several fundamental algebraic and analytic properties;
once these are proved, one can use it as a black box without caring too much for
the precise definition of $\Op_h(a)$. Of particular importance are
the Product, Adjoint, and Commutator Rules:
\begin{align}
\label{e:product-rule}
\Op_h(a)\Op_h(b)&=\Op_h(ab)+\mathcal O(h),\\
\label{e:adjoint-rule}
\Op_h(a)^*&=\Op_h(\bar a)+\mathcal O(h),\\
\label{e:commutator-rule}
[\Op_h(a),\Op_h(b)]&=-ih\Op_h(\{a,b\})+\mathcal O(h^2),
\end{align}
and the $L^2$ boundedness statement: if $a\in \CIc(T^*M)$ then
$\lVert\Op_h(a)\rVert_{L^2\to L^2}$ is bounded uniformly in~$h$.

\subsection{Semiclassical measures for eigenfunctions}
\label{s:semimes}

We can now introduce the main object of study in this article, which
are semiclassical measures associated to high frequency sequences
of eigenfunctions of the Laplacian. Semiclassical measures were originally introduced independently
by G\'erard~\cite{Gerard-measures} and Lions--Paul~\cite{Lions-Paul}.
We refer to~\cite[Chapter~5]{Zworski-Book} for a detailed treatment.

Following~\eqref{e:h-chosen},
we write the eigenvalue as $h^{-2}$ where $h$ is small.
Let $(M,g)$ be a Riemannian
manifold and consider a sequence of Laplacian eigenfunctions:
\[
-\Delta_g u_j =h_j^{-2} u_j,\qquad
h_j\to 0,\qquad
u_j\in C^\infty(M),\qquad
\lVert u_j\rVert_{L^2}=1.
\]
\begin{defi}
\label{d:weak-limit-2}
We say that the sequence $u_j$ converges semiclassically to a finite Borel measure
$\mu$ on the cotangent bundle $T^*M$, if
\begin{equation}
  \label{e:weak-limit-2}
\langle\Op_{h_j}(a)u_j,u_j\rangle_{L^2}\to \int_{T^*M}a(x,\xi)\,d\mu(x,\xi)\quad\text{as}\quad
j\to \infty
\end{equation}
for all test functions $a\in \CIc(T^*M)$. A measure $\mu$ on $T^*M$
is called a \emph{semiclassical measure} if it is the limit
of some sequence of Laplacian eigenfunctions.
\end{defi}
The statement~\eqref{e:weak-limit-2} actually applies to a broader class of symbols $a$
with polynomial growth as $\xi\to\infty$.
By~\eqref{e:op-h-mul}, if $a(x,\xi)=a(x)$ depends only on the position variable~$x$,
then the left-hand side of~\eqref{e:weak-limit-2}
is the integral $\int_M a|u_j|^2\,d\vol_g$. Comparing~\eqref{e:weak-limit-2}
with~\eqref{e:weak-limit-1}, we see that if $u_j$ converges semiclassically
to $\mu$, then it converges weakly to the pushforward of $\mu$ to the base $M$.
Thus we can think of semiclassical measures as (microlocal) lifts of the
weak limits of Definition~\ref{d:weak-limit-1}.

A quantum mechanical interpretation of semiclassical measures is as follows:
if $a\in C^\infty(T^*M)$ is a \emph{classical observable} (a function of position
and momentum) then $\Op_h(a)$ is the corresponding \emph{quantum observable}
and the expression $\langle\Op_h(a)u,u\rangle_{L^2}$
is the average value of the observable $a$ on the quantum particle with wave function~$u$.
Thus~\eqref{e:weak-limit-2} gives macroscopic information on the concentration of the particle
in both position and momentum in the high energy limit. Recalling~\eqref{e:op-h-fourier}, we can also interpret semiclassical measures as capturing the concentration of $u_j$
simultaneously in position and frequency.

One important property of Definition~\ref{d:weak-limit-2} is the presence
of compactness: any sequence of eigenfunctions has a subsequence
converging semiclassically to some measure~-- see~\cite[Theorem~5.2]{Zworski-Book}
and~\cite[Theorem~E.42]{DZ-Book}.
Other basic properties of semiclassical measures are summarized in the following
\begin{prop}
\label{l:measures-basic}
Let $\mu$ be a semiclassical measure for a Riemannian manifold $(M,g)$. Then:
\begin{itemize}
\item $\mu$ is a probability measure;
\item $\mu$ is supported on the cosphere bundle
\[
S^*M:=\{(x,\xi)\in T^*M\colon |\xi|_g=1\};
\]
\item $\mu$ is invariant under the geodesic flow
\[
\varphi^t:S^*M\to S^*M.
\]
Here the geodesic flow is naturally a flow on the sphere bundle $SM$,
which is identified with $S^*M$ using the metric~$g$.
\end{itemize}
\end{prop}
We give a sketch of the proof of Proposition~\ref{l:measures-basic}
to show how the fundamental properties~\eqref{e:product-rule}--\eqref{e:commutator-rule}
can be used. The first claim follows by taking $a=1$ in~\eqref{e:weak-limit-2},
in which case $\Op_h(a)$ is the identity operator. To see the second claim,
we use that the semiclassically rescaled Laplacian $-h^2\Delta_g$ is a quantization
of the quadratic function $|\xi|^2_g$ (giving the square of the length
of the cotangent vector~$\xi\in T_x^*M$ with respect to the metric~$g$), so
\[
P(h):=-h^2\Delta_g-1=\Op_h(|\xi|^2_g-1)+\mathcal O(h),\qquad
P(h_j)u_j=0.
\]
Now if $a\in \CIc(T^*M)$ vanishes on $S^*M$, we can write $a=b(|\xi|^2_g-1)$
for some $b\in \CIc(T^*M)$. By the Product Rule~\eqref{e:product-rule}
\[
\Op_{h_j}(a)u_j=\Op_{h_j}(b)P(h_j)u_j+\mathcal O(h_j)=\mathcal O(h_j)
\]
which by~\eqref{e:weak-limit-2} gives $\int_{T^*M}a\,d\mu=0$. Since
this is true for any $a$ vanishing on $S^*M$, we see that $\supp\mu\subset S^*M$ as needed.

The last claim is also simple to prove: if $b\in \CIc(T^*M)$ is arbitrary,
then
\[
0=\langle [P(h_j),\Op_{h_j}(b)]u_j,u_j\rangle_{L^2}
=-ih_j\langle \Op_{h_j}(\{|\xi|_g^2,b\})u_j,u_j\rangle_{L^2}+\mathcal O(h_j^2).
\]
Here the first equality follows from the fact that $P(h_j)u_j=0$
and $P(h_j)$ is self-adjoint; the second one uses the Commutator Rule~\eqref{e:commutator-rule}. Now~\eqref{e:weak-limit-2} shows that
the Poisson bracket $\{|\xi|_g^2,b\}$ integrates to~0 with respect to $\mu$.
But the Hamiltonian flow of $|\xi|_g^2/2$, restricted to $S^*M$, is 
the geodesic flow $\varphi^t$, so we get
\[
\int_{S^*M} \partial_t|_{t=0}(b\circ\varphi^t)\,d\mu=0\quad\text{for all}\quad
b\in \CIc(T^*M)
\]
from which it follows that $\int_{S^*M} b\circ\varphi^t\,d\mu$ is independent
of~$t$ and thus $\mu$ is invariant under the flow $\varphi^t$.

We now give the microlocal formulation of the question~\eqref{e:question-1}
asked at the beginning of the article:
\begin{equation}
\label{e:question-2}
\begin{gathered}
\text{For a given Riemannian manifold $(M,g)$, what can we say}\\
\text{about the set of all semiclassical measures?}
\end{gathered}
\end{equation}
The general expectation is that
\begin{itemize}
\item when the geodesic flow on $(M,g)$ is
`predictable', i.e. completely integrable, there are
semiclassical measures which can concentrate on small flow-invariant sets;
\item on the other hand, when the geodesic flow on $(M,g)$
has chaotic behavior, semiclassical measures
have to be more `spread out'.
\end{itemize}
One of the results supporting the first point above is the following
theorem of Jakobson--Zelditch~\cite{Jakobson-Zelditch}:
if $M$ is the round sphere then \emph{any} measure satisfying the conclusions
of Proposition~\ref{l:measures-basic} is a semiclassical measure.
See also the work of Studnia~\cite{Studnia-Harmonic-Oscillator} 
and Arnaiz--Maci\`a~\cite{Arnaiz-Macia} in the related
case of the quantum harmonic oscillator.

The rest of this article presents various results which support the
second point above, in particular giving several ways of defining
chaotic behavior of the geodesic flow and the way in which a measure is `spread out'.

\section{Ergodic systems}
\label{s:QE}

We first describe what happens under a `mildly chaotic' assumption on the geodesic flow
$\varphi^t:S^*M\to S^*M$, namely that it is \emph{ergodic} with respect to the Liouville measure. Here the Liouville measure $\mu_L=cd\vol_g(x)\,dS(\xi)$ is a natural flow-invariant probability measure on $S^*M$, with $dS$ denoting the volume measure on the sphere
$S_x^*M$ corresponding to~$g$ and $c$ some constant. By definition,
the flow $\varphi^t$ is ergodic with respect to $\mu_L$ if every $\varphi^t$-invariant Borel
subset $\Omega\subset S^*M$ has $\mu_L(\Omega)=0$ or $\mu_L(\Omega)=1$.

\begin{figure}
\includegraphics[width=15cm]{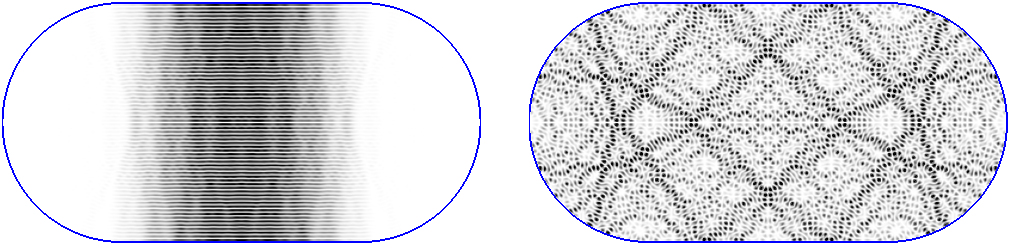}
\caption{Two Dirichlet eigenfunctions for a Bunimovich stadium, courtesy
of Alex Barnett (see the caption to Figure~\ref{f:billiards}): the right one shows equidistribution but the left
one does not. Quantum Ergodicity implies that most eigenfunctions look
from afar like the one on the right.}
\label{f:Bunimovich}
\end{figure}

We say that a sequence of eigenfunctions $u_j$ \emph{equidistributes}
if it converges to $\mu_L$ in the sense of Definition~\ref{d:weak-limit-2};
that is, in the high energy limit the probability of finding the corresponding
quantum particle in a set becomes proportional to the volume of this set.
A central result in quantum chaos is the following
Quantum Ergodicity theorem of Shnirelman~\cite{Shnirelman1},
Zelditch~\cite{Zelditch-QE}, and Colin de Verdi\`ere~\cite{CdV-QE},
which states that when the geodesic flow is ergodic, most eigenfunctions
equidistribute:
\begin{theo}
\label{t:QE}
Assume that the geodesic flow is ergodic with respect to the Liouville measure.
Then for any choice of orthonormal basis of eigenfunctions
$\{u_k\}$ there exists a density 1 subsequence $u_{k_j}$ which converges
semiclassically to~$\mu_L$ in the sense of Definition~\ref{d:weak-limit-2}.
\end{theo}
See~\cite[Chapter~15]{Zworski-Book} and the review of Dyatlov~\cite{around-shnirelman} for more recent expositions of the proof.
The version of Theorem~\ref{t:QE} for compact manifolds with boundary was proved
by G\'erard--Leichtnam~\cite{Gerard-Leichtnam} for convex domains in $\mathbb R^n$ with $W^{2,\infty}$ boundaries
and Zelditch--Zworski~\cite{Zelditch-Zworski} for compact Riemannian manifolds
with piecewise $C^\infty$ boundaries.
In this setting one imposes (Dirichlet or Neumann) boundary conditions on the eigenfunctions
and the geodesic flow is naturally replaced by the billiard ball flow (reflecting
off the boundary). See Figures~\ref{f:billiards} and~\ref{f:Bunimovich}
for numerical illustrations.

A natural question is whether the entire sequence of eigenfunctions
equidistributes, i.e. whether $\mu_L$ is the \emph{only} semiclassical
measure. For general manifolds with ergodic classical flows this is not always
true, as proved by Hassell~\cite{Hassell-Bunimovich}. In particular, for the
case of the Bunimovich stadium shown on Figure~\ref{f:Bunimovich}
the paper~\cite{Hassell-Bunimovich} shows that for almost every choice
of the parameter of the stadium (i.e. the aspect ratio of its central rectangle)
there exist semiclassical measures which are not the Liouville measure.

Another natural question is what happens when the classical flow has
\emph{mixed} behavior, e.g. $S^*M$ is the union of two
flow-invariant sets of positive Lebesgue measure such that the flow is ergodic on one of them and completely integrable on the other. \emph{Percival's Conjecture} claims that this mixed
behavior translates to macroscopic behavior of eigenfunctions, namely
one can split any orthonormal basis of eigenfunctions into three parts:
one of them equidistributes in the ergodic region,
another has semiclassical measures supported in the completely integrable region,
and the remaining part has density~0. A version of this conjecture for mushroom billiards
was proved by Gomes in his thesis~\cite{Gomes-thesis,Gomes-mushrooms}; see also the earlier work
of Galkowski~\cite{Galkowski-mushrooms} and Rivi\`ere~\cite{Riviere-mushrooms}.

\section{Strongly chaotic systems}
\label{s:anosov}

\begin{figure}
\includegraphics[height=6.5cm]{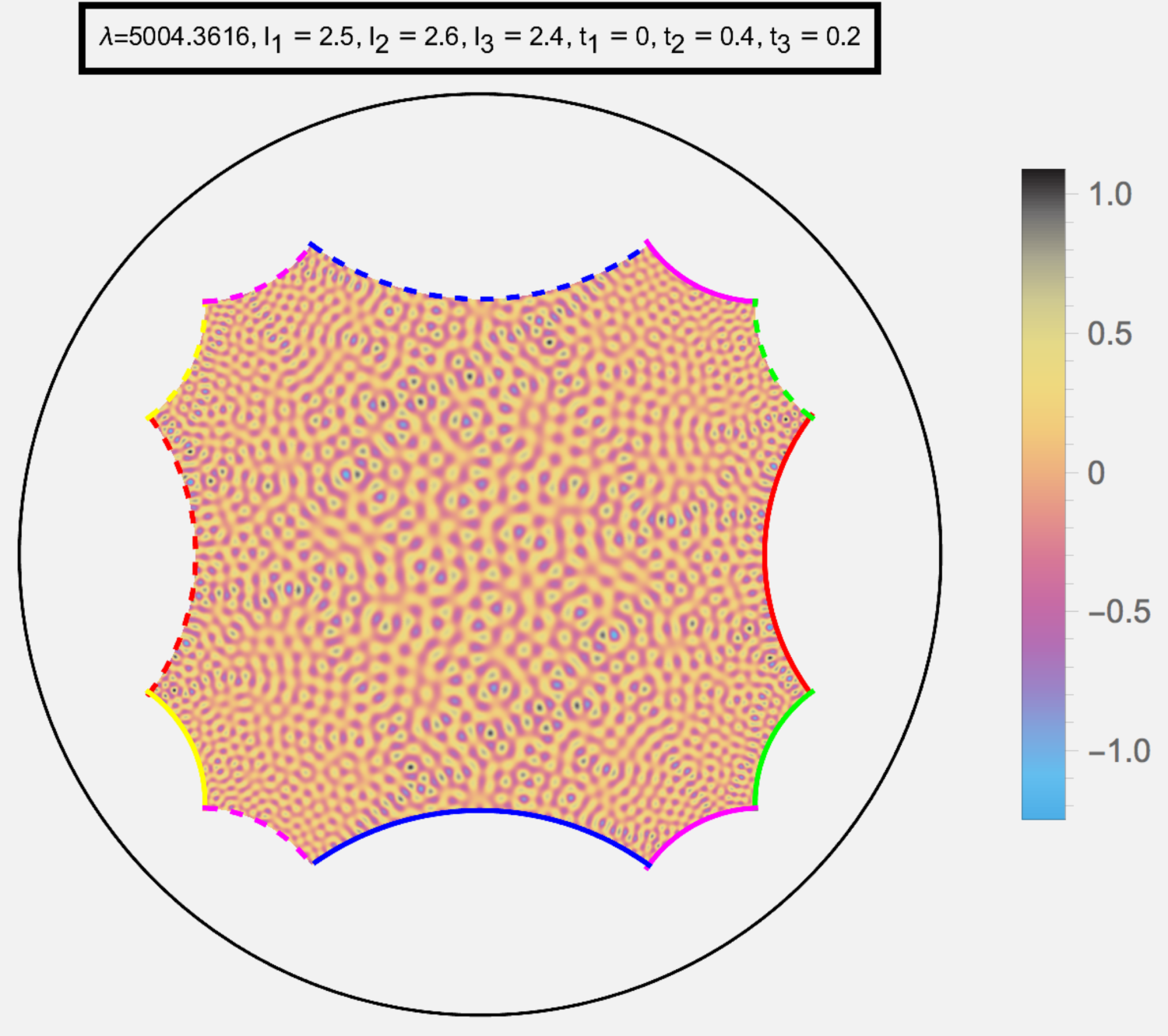}\quad
\includegraphics[height=6.5cm]{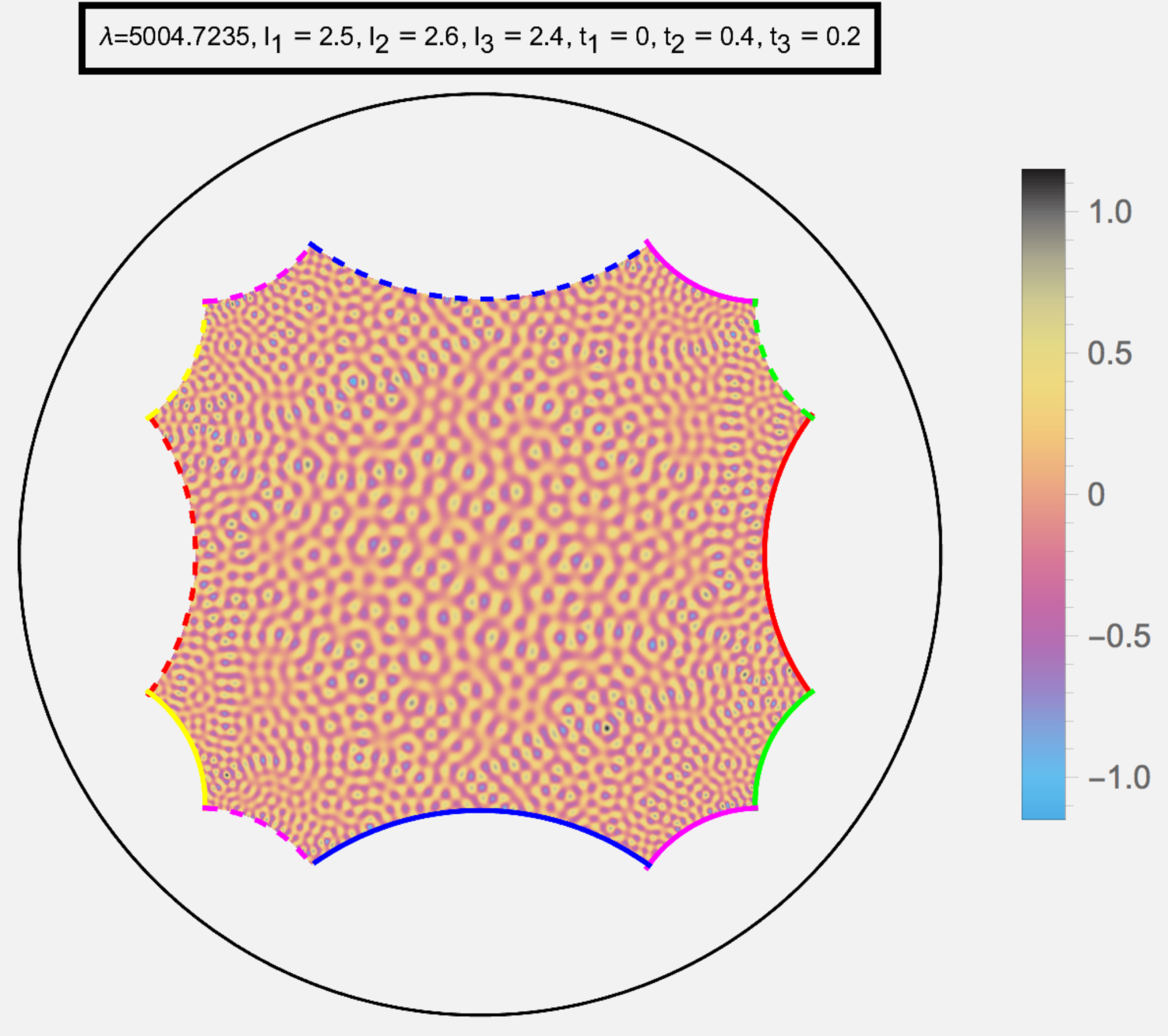}
\caption{Two Laplacian eigenfunctions on a hyperbolic surface, courtesy
of Alex Strohmaier (see Strohmaier--Uski~\cite{Strohmaier-Uski}). Here we view the surface
as a quotient of the hyperbolic plane by a group of isometries, or equivalently
as the result of gluing together appropriate sides of the pictured fundamental domain.
On a microscopic level the two eigenfunctions look different, but the macroscopic features are the same~-- both show equidistribution.}
\label{f:stro}
\end{figure}

We now describe what is known when the geodesic flow on $M$ is assumed to be
strongly chaotic. The latter assumption is understood in the sense of the following
\emph{Anosov property}:
\begin{defi}
Let $(M,g)$ be a compact Riemannian manifold without boundary.
We say that the geodesic flow $\varphi^t:S^*M\to S^*M$ has the Anosov property
if there exists a flow/unstable/stable decomposition
of the tangent spaces
\[
T_\rho (S^*M)=E_0(\rho)\oplus E_u(\rho)\oplus E_s(\rho),\quad
\rho\in S^*M,
\]
where $E_0$ is the one dimensional space spanned by the generator of the flow and
$E_u,E_s$ depend continuously on~$\rho$, are invariant under the flow $\varphi^t$,
and satisfy the exponential decay condition for some $\theta>0$:
\[
|d\varphi^t(\rho)v|\leq Ce^{-\theta|t|}|v|,\quad
\begin{cases}
v\in E_u(\rho),& t\leq 0;\\
v\in E_s(\rho),& t\geq 0.
\end{cases}
\]
\end{defi}
A large family of manifolds with Anosov geodesic
flows is given by compact Riemannian manifolds of negative sectional curvature,
see the book of Anosov~\cite{Anosov-book}. An important special case
is given by \emph{hyperbolic surfaces}, which
are compact oriented Riemannian manifolds of dimension~2
with Gauss curvature identically equal to~$-1$. See Figure~\ref{f:stro}
for a numerical illustration.

The Anosov property implies that the geodesic flow is ergodic with respect to the Liouville
measure, so Quantum Ergodicity applies to give that most eigenfunctions equidistribute.
The major open question is the following \emph{Quantum Unique Ergodicity} conjecture
which claims equidistribution for the entire sequence of eigenfunctions:
\begin{conj}
\label{c:QUE}
Assume that $(M,g)$ is a compact Riemannian manifold with Anosov geodesic flow.
Then $\mu_L$ is the only semiclassical measure.
\end{conj}
Conjecture~\ref{c:QUE} was originally stated by Rudnick--Sarnak~\cite{Rudnick-Sarnak-QUE}
for negatively curved Riemannian manifolds.
It is known in the special case of \emph{arithmetic} hyperbolic surfaces,
which have additional symmetries commuting with the Laplacian, called Hecke operators,
and we consider a joint basis of eigenfunctions of the Laplacian and a Hecke operator~--
see Lindenstrauss~\cite{Lindenstrauss-QUE} and Brooks--Lindenstrauss~\cite{Brooks-Lindenstrauss-QUE}. In general, in spite of significant partial progress described below,
the conjecture is open. One of the issues with a potential proof
is that Quantum Unique Ergodicity fails in the related setting of quantum cat maps~-- see Theorem~\ref{t:cat-que-fails} below.

\subsection{Entropy bounds}

A major step towards Quantum Unique Ergodicity (Conjecture~\ref{c:QUE})
are \emph{entropy bounds}, originating in the work of Anantharaman~\cite{Anantharaman-Entropy}:
\begin{theo}
\label{t:entropy-1}
Assume that the geodesic flow on $(M,g)$ has the Anosov property. Then
any semiclassical measure $\mu$ has positive Kolmogorov--Sinai entropy:
$\mathbf h_{\mathrm{KS}}(\mu)>0$.
\end{theo}
Here the Kolmogorov--Sinai entropy $\mathbf h_{\mathrm{KS}}(\mu)$ is a nonnegative number associated
to each flow-invariant measure~$\mu$; roughly speaking it expresses the
complexity of the flow from the point of view of that measure,
and is one way to measure how `spread out' the measure is~-- measures
which are more concentrated have lower entropy, and measures which
are more spread out have higher entropy.
Theorem~\ref{t:entropy-1} in particular implies the following
conjecture of Colin de Verdi\`ere~\cite{CdV-QE}:
\begin{equation}
  \label{e:QUE-conjecture}
\begin{gathered}
\text{On a hyperbolic surface, no semiclassical measure}\\
\text{can be supported on a closed geodesic}
\end{gathered}
\end{equation}
since the entropy of a measure supported on a closed geodesic is zero.

The lower bound on entropy in Theorem~\ref{t:entropy-1} is in general complicated.
However, in the case of hyperbolic (i.e. constant negative curvature) manifolds
Anantharaman--Nonnenmacher~\cite{Anantharaman-Nonnenmacher-Entropy} gave the following easy to state bound:
\begin{theo}
  \label{t:entropy-2}
Assume that $(M,g)$ is an $n$-dimensional hyperbolic manifold. Then
any semiclassical measure $\mu$ satisfies
\begin{equation}
  \label{e:entropy-2}
\mathbf h_{\mathrm{KS}}(\mu)\geq\textstyle{n-1\over 2}.
\end{equation}
\end{theo}
We remark that the Liouville measure in this setting has entropy $n-1$,
so~\eqref{e:entropy-2} in some sense excludes `half' of all invariant
measures as possible semiclassical measures.
For other entropy(-type) bounds, see the works of Anantharaman--Koch--Nonnenmacher~\cite{Anantharaman-Koch-Nonnenmacher}, Rivi\`ere~\cite{Riviere-Entropy-1,Riviere-Entropy-2}, and Anantharaman--Silberman~\cite{Anantharaman-Silberman}.

The constant in the bound~\eqref{e:entropy-2} matches (in the case of surfaces) the counterexamples for quantum cat maps
given in Theorem~\ref{t:cat-que-fails} below. Thus an important milestone on the way to
Quantum Unique Ergodicity would be to prove the following
\begin{conj}
Let $\mu$ be a semiclassical measure on an $n$-dimensional hyperbolic manifold~$(M,g)$.
Then $\mathbf h_{\mathrm{KS}}(\mu)>{n-1\over 2}$.
\end{conj}
We conclude this subsection with another conjecture which would go a long way towards
Quantum Unique Ergodicity but does not exclude the counterexample of Theorem~\ref{t:cat-que-fails}:
\begin{conj}
Let $\mu$ be a semiclassical measure on a compact manifold $(M,g)$
with Anosov geodesic flow. Then we have $\mu=\alpha\mu_L+(1-\alpha)\mu'$
for some $\alpha\in (0,1]$, where $\mu_L$ is the Liouville measure
and $\mu'$ is some probability measure on $S^*M$.
\end{conj}

\subsection{Full support property}
\label{s:full-support}

Another way to characterize how much a measure~$\mu$ is `spread out' is by looking
at its support, $\supp\mu\subset S^*M$. For surfaces with Anosov geodesic flows,
Dyatlov--Jin~\cite{meassupp} (in the hyperbolic case) and
Dyatlov--Jin--Nonnenmacher~\cite{varfup} (in the general case)
showed that the support of every semiclassical measure is the
entire $S^*M$:
\begin{theo}
  \label{t:meassupp}
Let $\mu$ be a semiclassical measure on a compact surface $(M,g)$ with 
Anosov geodesic flow. Then $\supp\mu=S^*M$, that is
$\mu(U)>0$ for every nonempty open set $U\subset S^*M$.
\end{theo}
Theorem~\ref{t:meassupp} and entropy bounds give different restrictions
on the set of possible semiclassical measures. On one hand (assuming $(M,g)$ is
a hyperbolic surface for simplicity), the entropy bound~\eqref{e:entropy-2}
implies that the Hausdorff dimension of $\supp\mu$ is at least~2,
but there exist flow-invariant measures supported on proper subsets
of $S^*M$ of dimension arbitrarily close to~3. On the other hand,
there exist measures which have full support and small entropy:
one can for example take a convex combination of the Liouville measure
and a measure supported on a closed geodesic.

The key new ingredient in the proof of Theorem~\ref{t:meassupp} is the
\emph{fractal uncertainty principle} of Bourgain--Dyatlov~\cite{fullgap}.
We state the following version appearing in~\cite{varfup}:
\begin{theo}
  \label{t:fup}
Let $\nu,h\in (0,1)$ and assume that $X,Y\subset\mathbb R$ are $\nu$-porous
up to scale~$h$, namely for any interval $I\subset\mathbb R$
of length $|I|\in [h,1]$ there exists a subinterval $J\subset I$
of length $|J|=\nu|I|$ such that $X\cap J=\emptyset$
(and similarly for $Y$). Then there exist constants
$C,\beta>0$ depending only on~$\nu$ such that
for all $f\in L^2(\mathbb R)$
\begin{equation}
  \label{e:fup}
\supp\hat f\subset h^{-1}Y\quad\Longrightarrow\quad
\lVert\mathbf 1_X f\rVert_{L^2(\mathbb R)}\leq Ch^\beta \lVert f\rVert_{L^2(\mathbb R)}.
\end{equation}
\end{theo}
One should think of the parameter $\nu$ in Theorem~\ref{t:fup} as fixed
and $h$ as going to~0. The sets $X,Y$ can depend on~$h$ as long as they
are $\nu$-porous; a basic example is given by $h\over 10$-neighborhoods
of some sets which are porous up to scale~0 (e.g. Cantor sets).
The estimate~\eqref{e:fup} can be interpreted as follows: if a function
$f$ lives in the (semiclassically rescaled) frequency space in a porous set~$Y$,
then only a small part of the $L^2$-mass of $f$ can concentrate on the
porous set~$X$. We refer the reader to the review~\cite{Dyatlov-JEDP}
for more details.

The proof of Theorem~\ref{t:meassupp} can be roughly summarized as follows
(restricting to the case of hyperbolic surfaces for simplicity):
assume that a sequence of eigenfunctions $\{u_j\}$ converges semiclassically to a
measure $\mu$ such that $\mu(\mathcal U)=0$ for some nonempty open set~$\mathcal U\subset S^*M$. Using microlocal methods, one can show that $u_j$
is in a certain sense concentrated on both of the sets
\[
\Omega_\pm(h_j):=\{\rho\in S^*M\mid \varphi^{\mp t}(\rho)\not\in \mathcal U\quad\text{for all}
\quad t\in [0,\log(1/h_j)]\}
\]
of geodesics which do not cross the set $\mathcal U$ in the future or in the past for
time $\log(1/h_j)$. Here one can barely make sense of localization in the position-frequency space
on each of the sets $\Omega_\pm(h_j)$, i.e. construct operators
$A_\pm$ which localize to these sets and write $u_j=A_+u_j+o(1)
=A_-u_j+o(1)$. However, the sets
$\Omega_\pm(h)$ have porous structure (see Figure~\ref{f:cat-holes} below for
the related case of quantum cat maps), and one can use the Fractal Uncertainty
Principle to show that $\lVert A_+A_-\rVert_{L^2\to L^2}=o(1)$, giving a contradiction.
We refer to~\cite{Dyatlov-JEDP} for a detailed exposition of the proof.

Theorem~\ref{t:meassupp} only applies to surfaces because the Fractal Uncertainty
Principle is only known for subsets of~$\mathbb R$. A na\"\i ve
generalization of Theorem~\ref{t:fup} to higher dimensions is false:
for example, the sets
\[
X=[0,h/10]\times [0,1],\ 
Y=[0,1]\times [0,h/10]\ \subset\ \mathbb R^2
\]
are both ${1\over 10}$-porous up to scale $h$ (where we replace intervals by
balls in the definition of porosity), but they do not
satisfy an estimate of type~\eqref{e:fup}: the Fourier transform of the
indicator function of $h^{-1}Y$ has large $L^2$ mass on~$X$.
(See~\cite[\S6]{FUP-ICMP} for a more detailed discussion.)
However, this does not translate to a counterexample for semiclassical measures,
leaving the door open for the following
\begin{conj}
\label{c:higher}
Let $\mu$ be a semiclassical measure on a compact manifold $(M,g)$
with Anosov geodesic flow. Then $\supp\mu=S^*M$.
\end{conj}
An analog of Conjecture~\ref{c:higher} is known for certain quantum cat maps,
see Theorem~\ref{t:highcat} below.

\section{Quantum cat maps}
\label{s:quantum-cat}

We finally discuss \emph{quantum cat maps}, which are toy models in quantum chaos
with microlocal properties similar to Laplacians on hyperbolic manifolds
(though the extensive research on them demonstrates that they are a `tough toy to crack'). They were
originally introduced by Hannay and Berry in~\cite{Hannay-Berry}.
We start
with two-dimensional quantum cat maps which are analogous to hyperbolic surfaces.
These maps quantize toral automorphisms (a.k.a. `Arnold cat maps')
\begin{equation}
  \label{e:cat-map}
x\mapsto Ax\bmod \mathbb Z^2,\quad
x\in \mathbb T^2=\mathbb R^2/\mathbb Z^2
\end{equation}
where $A\in\SL(2,\mathbb Z)$ is a $2\times 2$ integer matrix with determinant~1.
We make the assumption that $A$ is \emph{hyperbolic}, i.e. it has
no eigenvalues on the unit circle. A basic example of such a matrix is
\begin{equation}
  \label{e:basic-cat}
A=\begin{pmatrix} 2 & 1 \\ 1 & 1 \end{pmatrix}.
\end{equation}
Quantizations of the map~\eqref{e:cat-map} are not operators on $L^2$ of a manifold,
instead they are unitary $N\times N$ matrices, where the integer $N$ is related
to the semiclassical parameter $h$ as follows:
\[
2\pi N h = 1.
\]
The semiclassical limit $h\to 0$ studied above now turns
into the limit $N\to\infty$.

Before introducing quantizations of cat maps, we briefly discuss the adaptation of the quantization procedure~\eqref{e:quant-proc} to this setting, which has the form
\begin{equation}
  \label{e:Op-N}
a\in C^\infty(\mathbb T^2)\quad\mapsto\quad
\Op_N(a):\mathbb C^N\to\mathbb C^N.
\end{equation}
That is, functions on the 2-torus are quantized to $N\times N$ matrices.
The quantization procedure also depends on a twist parameter $\theta\in\mathbb T^2$,
but we suppress this in the notation. (If $N$ is even, then we can always just take
$\theta=0$ in what follows.)
See for example~\cite[\S2.2]{highcat} for more details.

Now, for $A\in\SL(2,\mathbb Z)$, its quantization is a family of unitary $N\times N$
matrices $B_N:\mathbb C^N\to\mathbb C^N$ which satisfies
the following \emph{exact Egorov's theorem}:
\begin{equation}
  \label{e:egorov-cat}
B_N^{-1}\Op_N(a)B_N=\Op_N(a\circ A)\quad\text{for all}\quad
a\in C^\infty(\mathbb T^2).
\end{equation}
Such $B_N$ exists and is unique modulo multiplication by a unit length scalar.
The statement~\eqref{e:egorov-cat} intertwines conjugation by $B_N$
(corresponding to quantum evolution) with pullback by the map~\eqref{e:cat-map}
(corresponding to classical evolution). It is analogous to Egorov's Theorem
for Riemannian manifolds (see e.g.~\cite[Theorem~15.2]{Zworski-Book}), which states that
\[
e^{-ith\Delta_g/2}\Op_h(a)e^{ith\Delta_g/2}=\Op_h(a\circ\varphi^t)+\mathcal O(h)
\]
where the geodesic flow $\varphi^t:S^*M\to S^*M$ is extended to $T^*M$ as the Hamiltonian
flow of $|\xi|_g^2/2$.
Thus the quantum cat map $B_N$ should be thought of as an analog of the Schr\"odinger propagator
$e^{ith\Delta_g/2}$, eigenfunctions of $B_N$ are analogous to Laplacian eigenfunctions, and the dynamics of the geodesic flow in this setting is replaced by the dynamics of the map~\eqref{e:cat-map}.

Using the quantization~\eqref{e:Op-N}, we can define similarly to~\eqref{e:weak-limit-2} semiclassical measures
associated to sequences of eigenfunctions
\[
B_{N_j} u_j=\lambda_ju_j,\qquad
u_j\in\mathbb C^{N_j},\qquad
\lVert u_j\rVert_{\ell^2}=1,\qquad
N_j\to\infty.
\]
These are probability measures on $\mathbb T^2$ which are invariant under the map~\eqref{e:cat-map} (as can be seen directly from Egorov's theorem~\eqref{e:egorov-cat}).

When the matrix $A$ is hyperbolic, the map~\eqref{e:basic-cat} is ergodic with respect to
the Lebesgue measure on $\mathbb T^2$. Using this fact, Bouzouina--de Bi\`evre~\cite{Bouzouina-deBievre} showed
Quantum Ergodicity in this setting: if we put together orthonormal bases
of eigenfunctions of $B_N$ for all~$N$, then there exists a density~1 subsequence
of this sequence which converges to the Lebesgue measure.

\begin{figure}
\includegraphics[height=7cm]{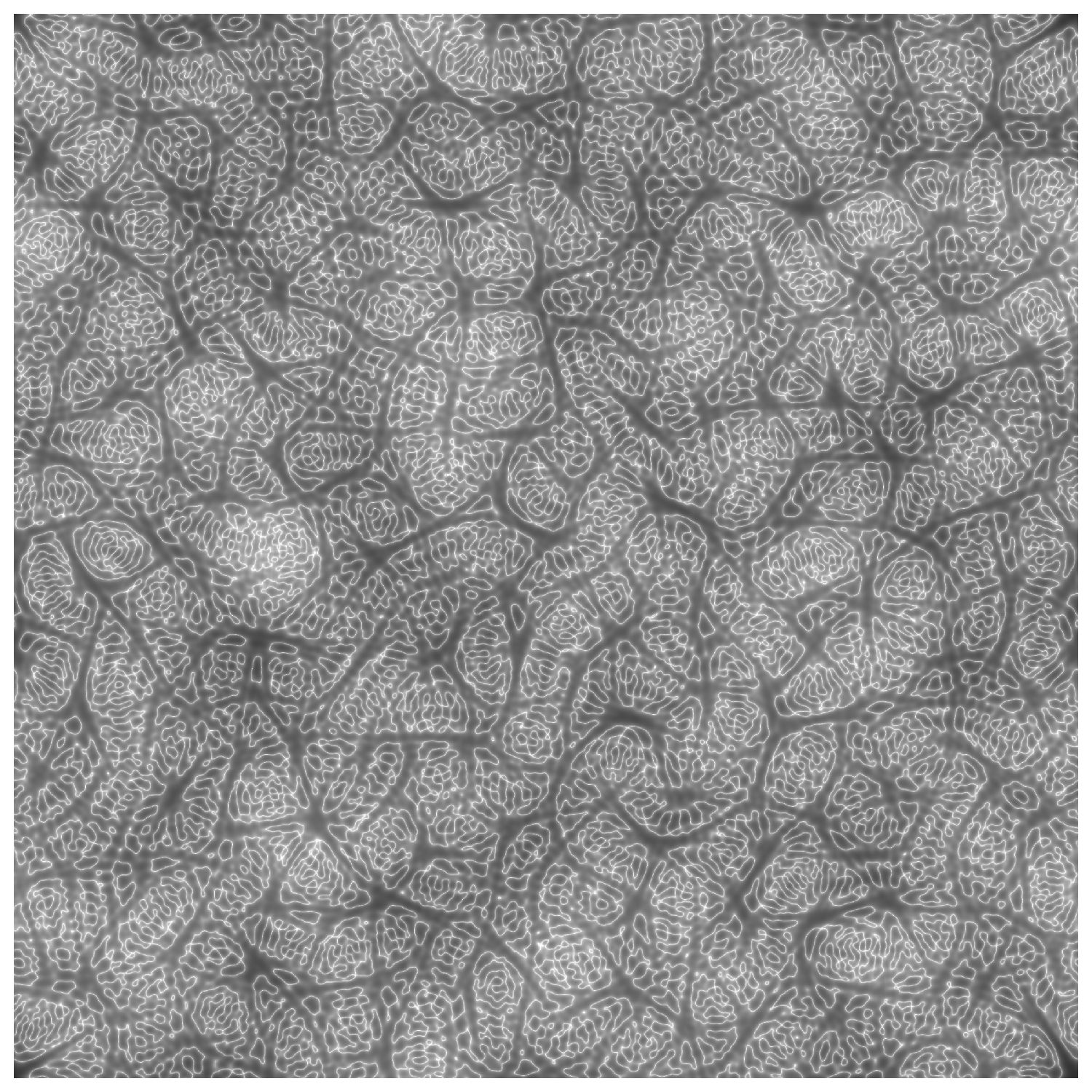}
\quad
\includegraphics[height=7cm]{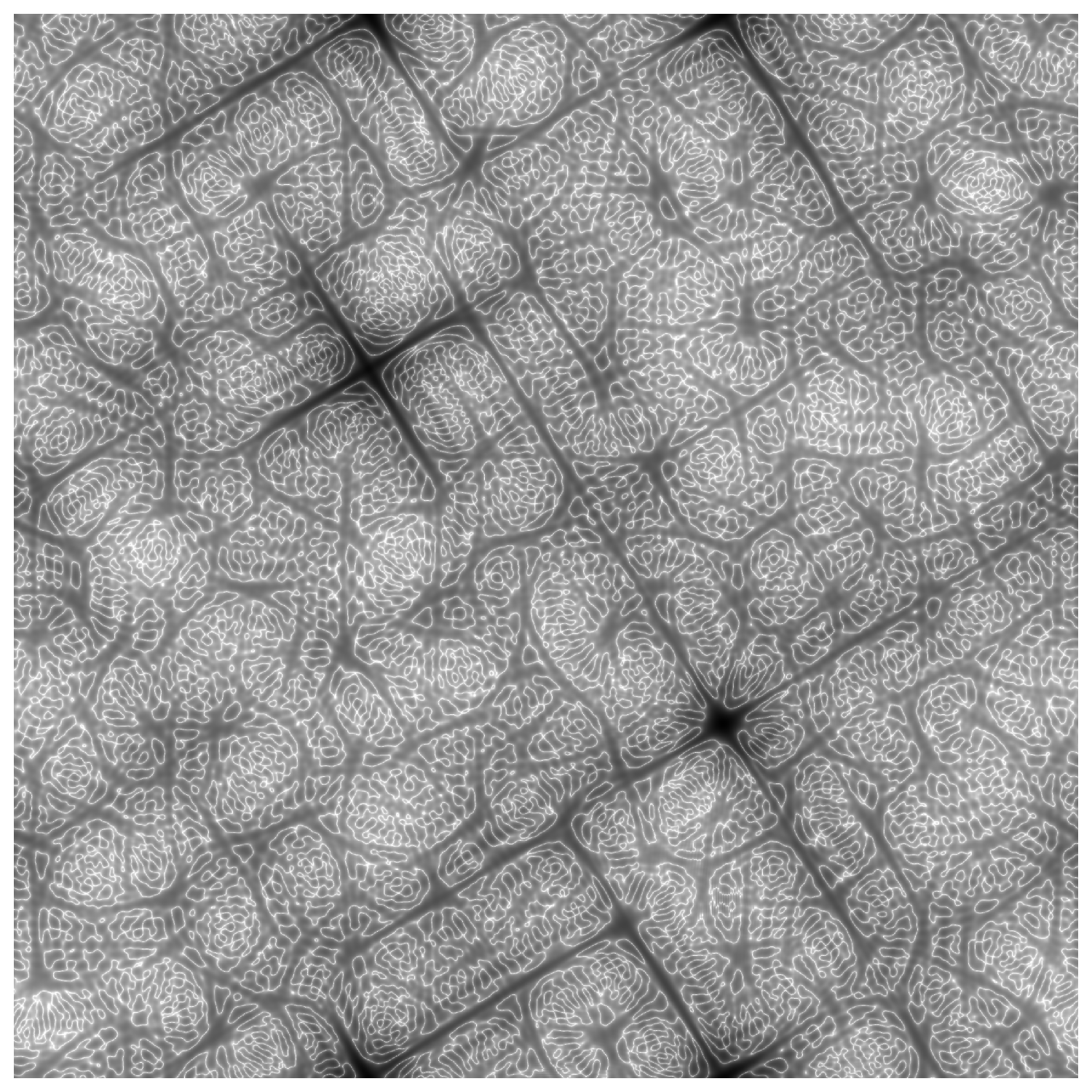}
\caption{Phase space concentration
for two eigenfunctions of the quantum cat map
with $A$ given by~\eqref{e:basic-cat}
and $N=1292$. More specifically,
we plot the absolute value of a smoothened out Wigner transform
of the eigenfunction on the logarithmic scale (see e.g.~\cite[\S2.2.5]{highcat}). On the left is a typical eigenfunction,
showing equidistribution. On the right is a particular
eigenfunction of the type constructed in~\cite{Faure-Nonnenmacher-dB},
corresponding to a measure of the type~\eqref{e:cat-que-fails}
featuring the closed trajectory $\{({1\over 3},0),({2\over 3},{1\over 3}),({2\over 3},0),({1\over 3},{2\over 3})\}$.
The existence of such an eigenfunction relies on the careful choice of~$N$:
$A^{18}$ is the identity matrix modulo $2N$.
}
\label{f:cat-eigs}
\end{figure}

On the other hand, Faure--Nonnenmacher--De Bi\`evre~\cite{Faure-Nonnenmacher-dB} showed that
Quantum Unique Ergodicity fails for quantum cat maps:
\begin{theo}
\label{t:cat-que-fails}
Let $A\in \SL(2,\mathbb Z)$ be a hyperbolic matrix. Fix
any periodic trajectory $\gamma\subset \mathbb T^2$ of the map~\eqref{e:cat-map}.
Then there exists
a sequence of eigenfunctions $u_j$ of the quantum cat map $B_{N_j}$,
for some $N_j\to\infty$, which converge semiclassically to the measure
\begin{equation}
  \label{e:cat-que-fails}
\textstyle{1\over 2}\delta_\gamma+\textstyle{1\over 2}\mu_L
\end{equation}
where $\delta_\gamma$ is the delta probability measure on the trajectory
$\gamma$ and $\mu_L$ is the Lebesgue measure on $\mathbb T^2$.
\end{theo}
We remark that the choice of $N_j$ in Theorem~\ref{t:cat-que-fails} is highly
special: one takes them so that the matrix $A^{k_j}$ is the identity modulo $2N_j$
where $k_j$ is very small, namely $k_j\sim \log N_j$. This implies
that the quantum cat map $B_{N_j}$ also has a short period,
namely $B_{N_j}^{k_j}$ is a scalar. See the papers of Dyson--Falk~\cite{Dyson-Falk}
and Bonechi--De Bi\`evre~\cite{Bonechi-deBievre} for more information on the periods
of the cat map. A numerical illustration of Theorem~\ref{t:cat-que-fails} is given on Figure~\ref{f:cat-eigs}.

The entropy of the measure~\eqref{e:cat-que-fails} is equal to
half the entropy of the Lebesgue measure. This matches
the constant in the entropy bound of Theorem~\ref{t:entropy-2}.
Since from the point of view of microlocal analysis quantum cat maps
have similar properties to hyperbolic surfaces, significant new insights would
be needed to show that a counterexample of the kind~\eqref{e:cat-que-fails}
cannot occur for hyperbolic surfaces.

Faure--Nonnenmacher~\cite{Faure-Nonnenmacher-cat} showed that the constant~${1\over 2}$ in~\eqref{e:cat-que-fails}
is sharp: the mass of the pure point part of any semiclassical measure for a quantum cat map is less than or equal to the mass
of its Lebesgue part. Brooks~\cite{Brooks-cat} generalized this to a statement that 
the mass of lower entropy components of any semiclassical measure is less than or equal to the mass of higher entropy
components; this in particular implies an entropy bound analogous to~\eqref{e:entropy-2}.

There is also an analogue of arithmetic Quantum Unique Ergodicity in the setting
of cat maps: Kurlberg--Rudnick~\cite{Kurlberg-Rudnick} introduced Hecke operators
which commute with $B_N$ and showed that any sequence of joint
eigenfunctions of $B_N$ and these operators converges to the Lebesgue measure.
This does not contradict the counterexample of Theorem~\ref{t:cat-que-fails}
since for the values of $N_j$ chosen there, the map $B_{N_j}$ has eigenvalues
of high multiplicity.

\begin{figure}
\includegraphics[height=4.75cm]{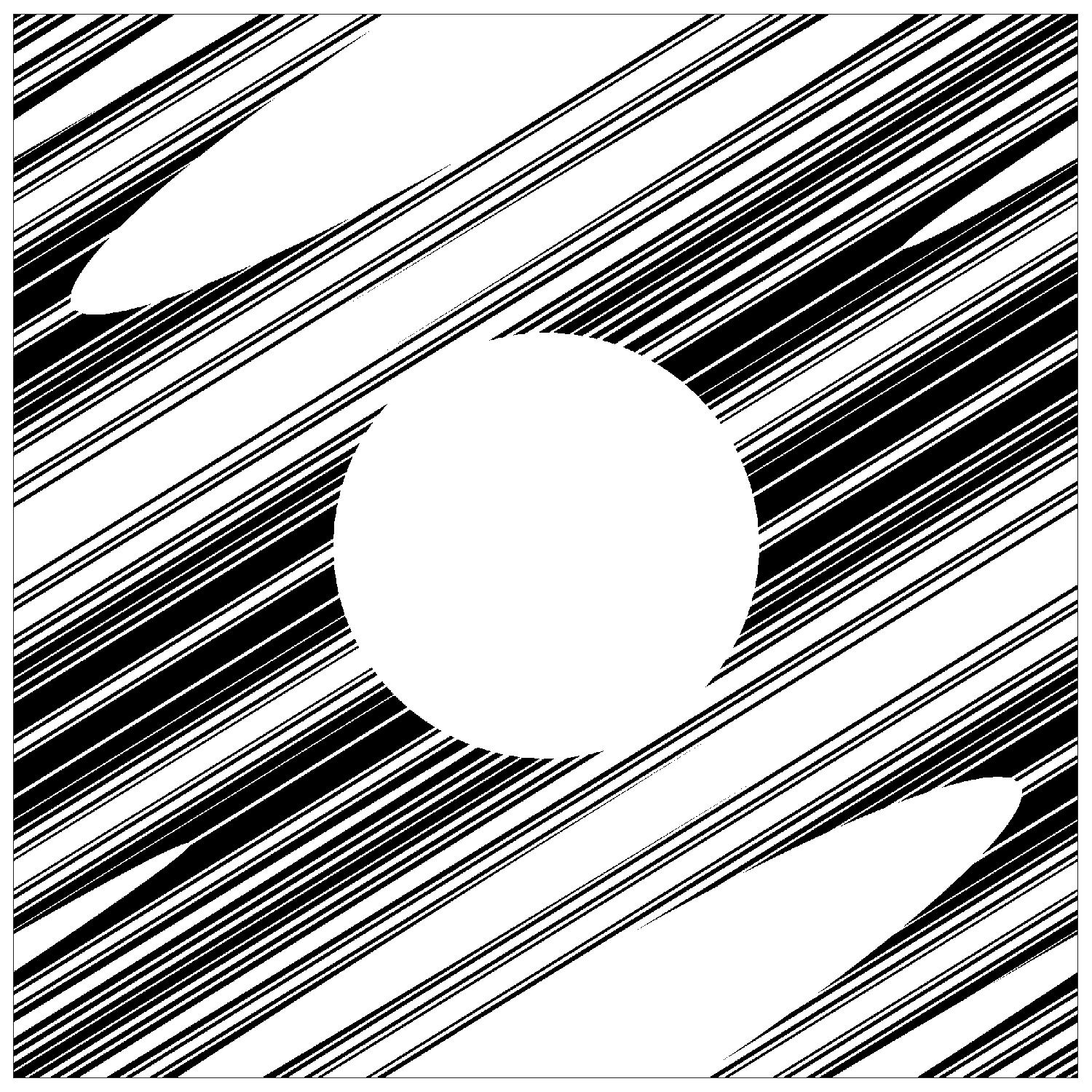}\quad
\includegraphics[height=4.75cm]{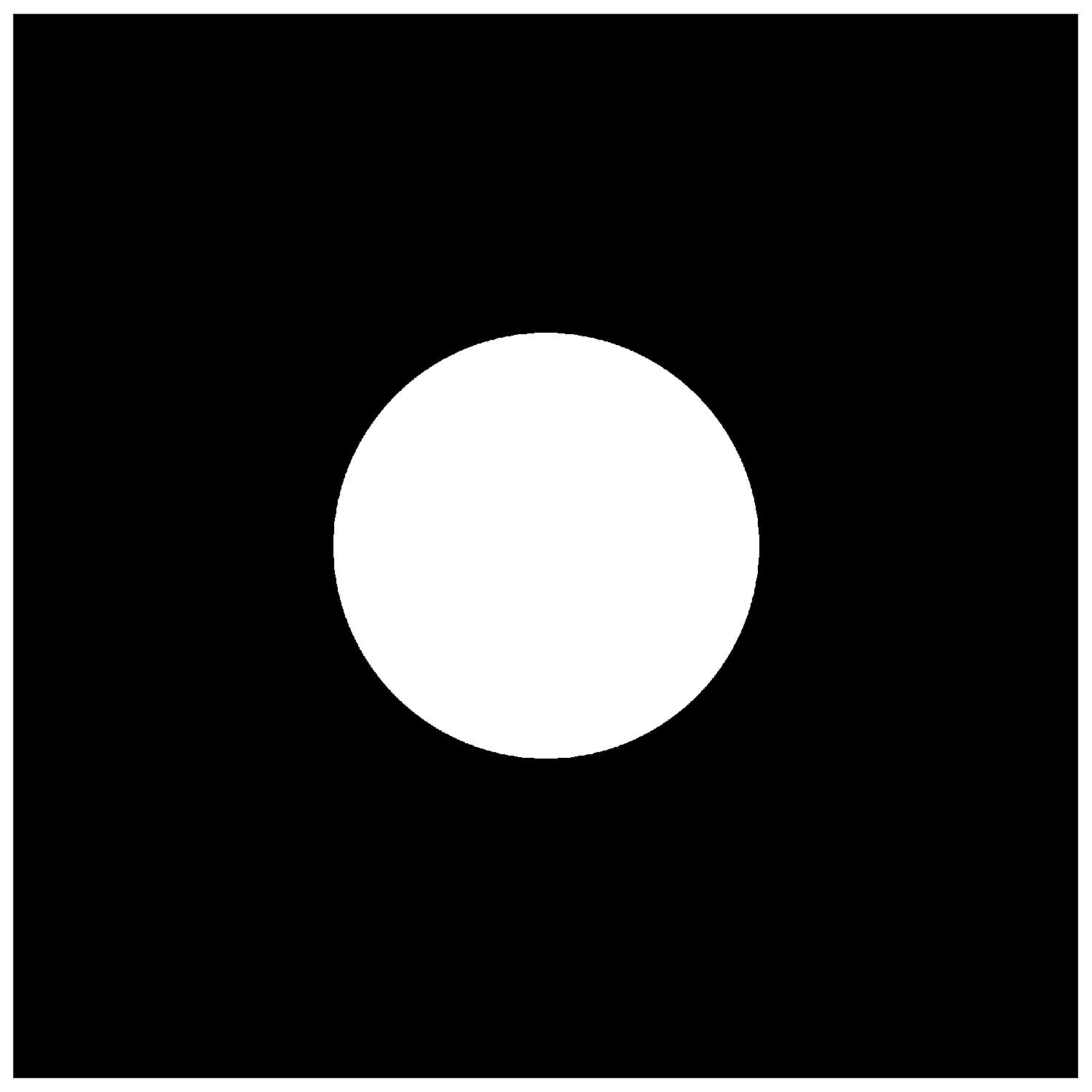}\quad
\includegraphics[height=4.75cm]{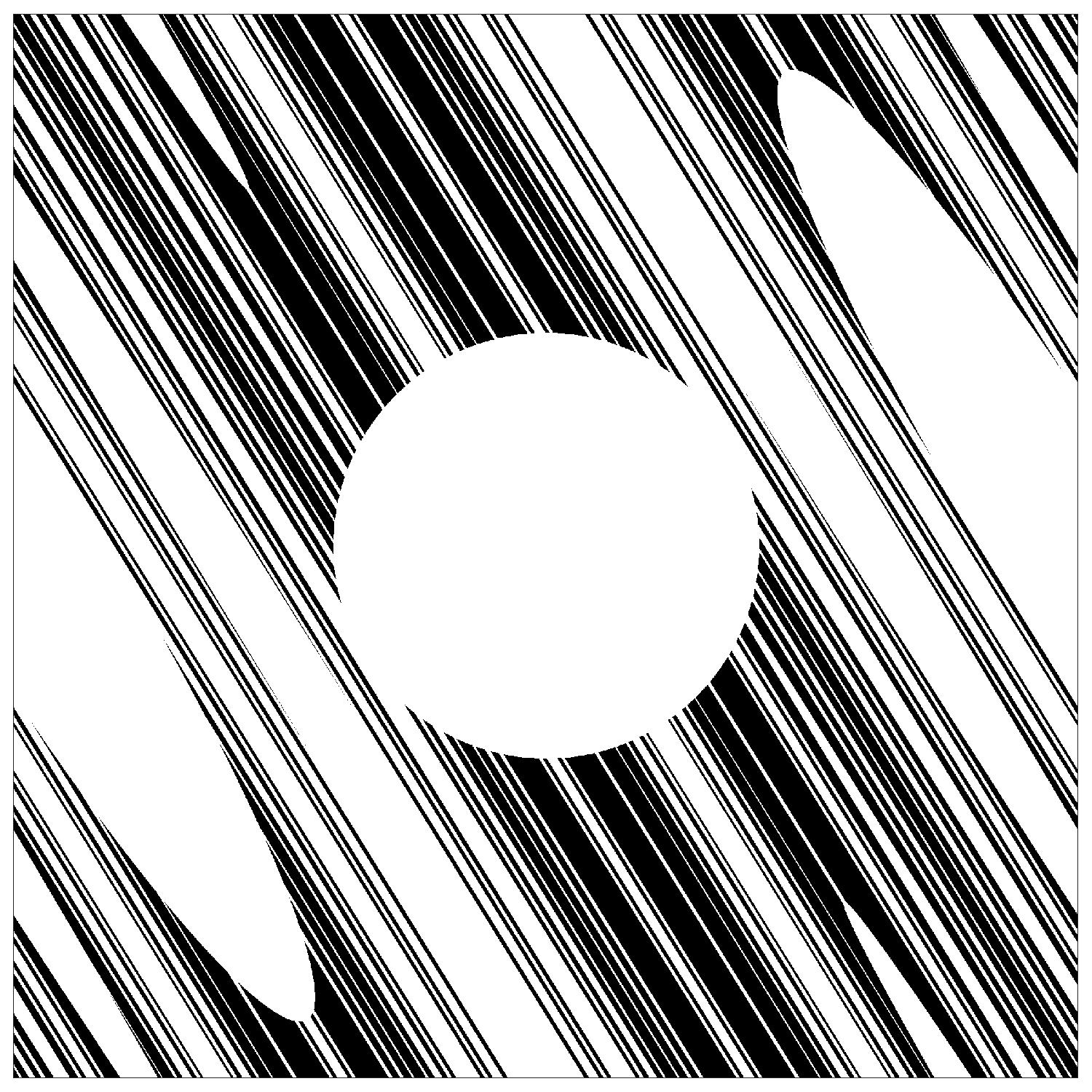}
\caption{A set $\mathcal U\subset\mathbb T^2$ (center picture, in white) and the corresponding
sets $\Omega_+(N),\Omega_-(N)$ (left/right picture). The set $\Omega_+(N)$ is `smooth' in the unstable direction
of the matrix $A$ and porous in the stable direction, with the porosity constant depending only on~$\mathcal U$. Same
is true for $\Omega_-(N)$ but switching the roles of the stable/unstable directions. The fractal uncertainty principle
of Theorem~\ref{t:fup} can be used to show that no function can be localized on both $\Omega_+(N)$ and $\Omega_-(N)$.}
\label{f:cat-holes}
\end{figure}

We now discuss the recent results on support of semiclassical
measures for cat maps, proved using the fractal uncertainty principle.
For two-dimensional cat maps, Schwartz~\cite{Schwartz-cat} showed the following
\begin{theo}
  \label{t:lowcat}
Let $\mu$ be a semiclassical measure for a quantum cat map associated
to some hyperbolic matrix $A\in\SL(2,\mathbb Z)$. Then
$\supp \mu=\mathbb T^2$.
\end{theo}
Similarly to~\S\ref{s:full-support}, the proof uses that no function can be localized simultaneously on
the two sets
\[
\Omega_\pm(N):=\bigg\{\rho\in \mathbb T^2\,\bigg|\, A^{\mp j}(\rho)\not\in\mathcal U\quad\text{for all}\quad
j=0,\dots,{\log N\over\log |\lambda_+|}\bigg\}
\]
where $\lambda_+$ is the eigenvalue of $A$ such that $|\lambda_+|>1$. Here $\mathcal U\subset\mathbb T^2$ is some nonempty open set. See Figure~\ref{f:cat-holes}.

We finally discuss the quantum cat map analog of the higher-dimensional Conjecture~\ref{c:higher},
by considering quantum cat maps associated to symplectic integer
matrices $A\in\Sp(2n,\mathbb Z)$. In this setting
Dyatlov--J\'ez\'equel~\cite{highcat} proved
\begin{theo}
  \label{t:highcat}
Let $\mu$ be a semiclassical measure for a quantum cat map associated
to a matrix $A\in\Sp(2n,\mathbb Z)$ such that:
\begin{itemize}
\item $A$ has a simple eigenvalue $\lambda_+$ such that
all other eigenvalues satisfy $|\lambda|<\lambda_+$; and
\item the characteristic polynomial of~$A$ is irreducible over the rationals.
\end{itemize}
Then $\supp\mu=\mathbb T^{2n}$.
\end{theo}
Here the first condition makes it possible to still use the one-dimensional fractal uncertainty principle in the proof.

We remark that there are examples of semiclassical measures
which do not have full support
for some matrices $A$ satisfying the first condition of Theorem~\ref{t:highcat} but not the second condition. In particular, there
exist semiclassical measures supported on tori associated
to any $A$-invariant rational Lagrangian subspace of $\mathbb R^{2n}$.
See the work of Kelmer~\cite{Kelmer-cat} and the discussion in~\cite[Appendix~A]{highcat}.


\medskip\noindent\textbf{Acknowledgements.}
The author was supported by NSF CAREER grant DMS-1749858
and a Sloan Research Fellowship.


\bibliographystyle{alpha}
\bibliography{General,Dyatlov,QC}

\end{document}